\documentclass[10pt,leqno]{article}
\usepackage{amsmath}

\usepackage{amssymb}

\usepackage{xypic}

\usepackage[all]{xy}

\usepackage{latexsym}

\usepackage{amsfonts}

\usepackage{epsf}

\usepackage[dvips]{graphicx}

\usepackage[latin1]{inputenc}

\usepackage[T1]{fontenc}

\newtheorem{teo}{Theorem}[section]

\newtheorem{lem}[teo]{Lemma}
\newtheorem{prop}[teo]{Proposition}

\newtheorem{df}[teo]{Definition}
\newtheorem{es}[teo]{Example}%%%
\newtheorem{oss}[teo]{Remark}%%%
\newcommand{\A}{\mathcal{A}}

\newcommand{\C}{\mathcal{C}}

\newcommand{\U}{\mathcal{U}}

\newcommand{\psh}{\mathrm{Psh}}

\newcommand{\ob}{\mathrm{Ob}}

\newcommand{\Ho}{\mathrm{Hom}}

\newcommand{\Cc}{\mathcal{C}^\land}

\newcommand{\iso}{\stackrel{\sim}{\to}}

\newcommand{\id}{\mathrm{id}}

\renewcommand{\dim}{\textbf{Proof.}}

\newcommand{\qed}{\nopagebreak \phantom{} \hfill $\Box$ \\}

\newcommand{\imin}[1]{#1^{-1}}

\newcommand{\lind}[1]{\underset{#1}{\underrightarrow{\lim}}}

\newcommand{\Lind}{\underrightarrow{\lim}}  %grazie Anna

\newcommand{\lpro}[1]{\underset{#1}{\underleftarrow{\lim}}}

\title{Proper stacks}
\author{Luca \textsc{Prelli}
          \footnote{Universit{\`a} di Padova, Dipartimento di Matematica Pura ed
Applicata, via Trieste 63, 35121 Padova, Italy or: Universit{\'e}
Pierre et Marie Curie, Institut de Math{\'e}matiques de Jussieu,
175 rue du Chevaleret, 75013 Paris, France e-mail address:
lprelli@math.unipd.it }}
\date{}

\begin{document}

\maketitle

\thispagestyle{empty}

\begin{abstract}
We generalize the notion of proper stack introduced by Kashiwara
and Schapira to the case of a general site, and we prove that a
proper stack is a stack.
\end{abstract}

\tableofcontents

\addcontentsline{toc}{section}{Introduction}

\section*{Introduction}

In \cite{KS01} Kashiwara and Schapira defined the notion of proper
stack on a locally compact topological space $X$. A proper stack
is a separated prestack $\mathcal{S}$ satisfying suitable
hypothesis. They proved that a proper stack is a stack. In this
paper, we generalize the notion of proper stack to the case of a
site $X$ associated to a small category $\C_X$ and we prove that a
proper
stack is a stack.\\

\noindent \textbf{Acknowledgments.} We thank Pierre Schapira who
encouraged us to generalize the notion of proper stack. We thank
St\'ephane Guillermou and Pietro Polesello for their useful
remarks.

\section{Review on Grothendieck topologies and sheaves}

 Let $\C$ be a category\footnote{We shall work in a given universe
 $\U$, small means $\U$-small (i.e. a set is $\U$-small if it is
isomorphic to a set belonging to $\U$) and a category $\C$ means a
$\U$-category (i.e. $\Ho_\C(X,Y)$ is $\U$-small for any $X,Y \in
\C$).}. As usual we denote by $\Cc$ the category of functors from
$\C^{op}$ to $\bf{Set}$ and we identify $\C$ with its image in
$\Cc$ via the Yoneda embedding. If $A \in \C^\land$, we will
denote by $\C_A$ the category of arrows $U \to A$ with $U \in \C$.
When taking inductive and projective limits on a category $I$ we
will always assume that $I$ is small.

We recall here some classical definitions (see \cite{SGA4}),
following the presentation of \cite{KS}.

%\begin{df} A presite $X$ is a category $\C_X$, and a morphism of presites $f:X \to Y$ is a functor $f^t:\C_Y \to \C_X$
%\end{df}

\begin{df}\label{LE} A Grothendieck topology on a small category $\C_X$ is a collection of morphisms in $\Cc_X$ called local epimorphisms, satysfying the following conditions:
\begin{itemize}
\item[LE1] For any $U \in \C_X$, $\id_U:U \to U$ is a local epimorphism.
\item[LE2] Let $A_1 \stackrel{u}{\to} A_2 \stackrel{v}{\to} A_3$ be morphisms in $\Cc_X$. If $u$ and $v$ are local epimorphisms, then $v \circ u$ is a local epimorphism.
\item[LE3] Let $A_1 \stackrel{u}{\to} A_2 \stackrel{v}{\to} A_3$ be morphisms in $\Cc_X$. If $v \circ u$ is a local epimorphism, then $v$ is a local epimorphism.
\item[LE4] A morphism  $u:A \to B$ in $\Cc_X$ is a local epimorphism if and only if for any $U \in \C_X$ and any morphism $U \to B$, the morphism $A \times_B U \to U$ is a local epimorphism.
\end{itemize}
\end{df}

\begin{df}
A morphism $A \to B$ in $\Cc_X$ is a local monomorphism if $A \to
A \times_B A$ is a local epimorphism.

A morphism $A \to B$ in $\Cc_X$ is a local isomorphism if it is
both a local epimorphism and a local monomorphism.
\end{df}

\begin{df} A site $X$ is a category $\C_X$ endowed with a Grothendieck
topology.
\end{df}

Let $\A$ be a category admitting small inductive and projective
limits.

\begin{df} An $\A$-valued presheaf on $X$ is a functor
$\C_X^{op} \to \A$. A morphism of presheaves is a morphism of such
functors. One denotes by $\psh(X,\A)$ the category of $\A$-valued
presheaves on $X$.
\end{df}

If $F \in \psh(X,\A)$, it extends naturally to $\Cc_X$ by setting
$$F(A)=\lpro {(U \to A) \in \C_A} F(U),$$
where $A \in \Cc_X$ and $U \in \C_X$.

\begin{df} Let $X$ be a site.
\begin{itemize}
\item One says that $F \in \psh(X,\A)$ is separated, if for any local isomorphism $A \to U$ with $U \in \C_X$ and $A \in \Cc_X$, $F(U) \to F(A)$ is a monomorphism.
\item  One says that $F \in \psh(X,\A)$ is a sheaf, if for any local isomorphism $A \to U$ with $U \in \C_X$ and $A \in \Cc_X$, $F(U) \to F(A)$ is an isomorphism.
%\item One denotes by $\sh(X,\A)$ the category of sheaves on $X$.
\end{itemize}
\end{df}

\section{Review on stacks}

Let $\C_X$ be a small category.
%If $A \in \Cc_X$ we denote by $\C_A$ the category of arrows $U \to A$ with $U \in \C_X$.
We suppose that a Grothendieck topology on $\C_X$ is defined and
we denote by $X$ the associated site.
%If $U_i,U_j,U \in \C_X$, with $U_i \to U$ and $U_j \to U$ we set for short $U_{ij}:=U_i\times_U U_j$.
We recall some classical definitions (see \cite{Gi71}), following
the presentation of \cite{KS}.
%\section{Definitions}

\begin{df} A prestack $\mathcal{S}$ on $X$ is the data of:

\begin{itemize}

\item for each $U\in\C_X$, a category $\mathcal{S}(U)$,

\item for each ${V}\to{U} \in \C_U$, a functor $j_{VU*}:{\mathcal{S}(U)}\to{\mathcal{S}(V)}$,

\item given $U,V,W \in \C_X$ and ${W}\to{V}\to{U}$, an
  isomorphism of functors $\lambda_{WVU}:{{j_{WV*}}\circ{j_{VU*}}}\stackrel{\sim}{\to}{j_{WU*}}$,

\end{itemize}

such that

\begin{itemize}

\item $j_{UU*}=\id_{\mathcal{S}(U)}$,

\item given $\{U_i\}_{i \in I} \in \C_X$, i=1,2,3,4 and
  ${U_1}\to{U_2}\to{U_3}\to{U_4}$, the following diagram commutes:

$$\xymatrix{{{j_{12*}}\circ{j_{23*}}\circ{j_{34*}}}\ar[rr]^{\lambda_{234}}\ar[d]^{\lambda_{123}}&&{{j_{12*}}\circ{j_{24*}}}\ar[d]^{\lambda_{124}}\\
{{j_{13*}}\circ{j_{34*}}}\ar[rr]^{\lambda_{134}}&&{j_{14*}}}$$

\end{itemize}

\end{df}

Let $\lpro {U \in \C_X}\mathcal{S}(U)$ denote a category defined
as follows. An object $F$ of $\lpro {U \in \C_X}\mathcal{S}(U)$ is
a family $\{(F_U)_U, (\psi_u)_u\}$ where
\begin{itemize}
\item for any $U \in \C_X$, $F_U \in \ob(\mathcal{S}(U))$,
\item for any morphism $U_1 \to U_2$ in $\C_X$, $\psi_{12}:j_{12*}F_{U_2}\to F_{U_1}$ is an isomorphism, such that for any sequence $U_1 \underset{}{\to} U_2 \underset{}{\to} U_3$ the following diagram
commutes

$$\xymatrix{j_{12*}j_{23*}F_{U_3} \ar[r]^{\psi_{23}} \ar[d]_{\lambda_{123}} & j_{12*}F_{U_2} \ar[d]^{\psi_{12}}\\
j_{13*}F_{U_3} \ar[r]^{\psi_{13}} & F_{U_1}.}$$ Note that
$\psi_{\id_U}=\id_{F_U}$ for any $U \in \C_X$.

\end{itemize}

\noindent The morphisms are defined in natural way. Let $F,G \in
\lpro {U \in \C_X}\mathcal{S}(U)$.
%Let $F=\{(F_U)_U, (\psi_u)_u\}$, $G=\{(G_U)_U, (\varphi_u)_u\}$.
Then
$$\Ho_{\lpro {U \in \C_X}\mathcal{S}(U)}(F,G)
\simeq \lpro {U \in \C_X} \Ho_{\mathcal{S}(U)}(F_U,G_U).$$

For any $A \in \Cc_X$, we set
$$\mathcal{S}(A)=\lpro {(U \to A) \in \C_A} \mathcal{S}(U)$$
%Then, if $X$ denotes the terminal object of $\Cc_X$, then $\Lpro \mathcal{S}=\mathcal{S}(X)$.
A morphism $\varphi:A \to B$ in $\Cc_X$ defines a functor
$j_{AB*}:\mathcal{S}(B) \to \mathcal{S}(A)$, therefore a prestack
on $\C_X$ extends naturally to a prestack on $\Cc_X$.

%We write for short $\imin {i_U}$ instead of $\imin {i_{UX}}$, and
%We will often write $F|_V$ instead of $\imin {i_{VU}}F$.

%\begin{df}  Let $U \in \C_X$, and $F,G \in \mathcal{S}(U)$. One denotes by
%$\ho_{\mathcal{S}|_U}(F,G)$ the presheaf $V \mapsto
%\Ho_{\mathcal{S}(V)}(j_{VU*}F,j_{VU*}G)$ for $V \to U$ in $\C_X$.
%\end{df}
%Note that for $A \in \Cc_X$ and $F,G \in \mathcal{S}(A)$, we obtain
%$$\Gamma(A;\ho_{\mathcal{S}|_A}(F,G)) \simeq \lpro {(U \to A) \in
%\C_A}\Ho_{\mathcal{S}(U)}(j_{UA*}F,j_{UA*}G).$$

\begin{df} Let $X$ be a site.
\begin{itemize}
%\item A prestack $\mathcal{S}$ on $X$ is called separated if for any $U \in \C_X$ and any $F,G \in \mathcal{S}(U)$, $\ho_{\mathcal{S}|_U}(F,G)$ is a sheaf on $U$.
\item A prestack $\mathcal{S}$ on $X$ is called separated if for any $U \in \C_X$, and for any local isomorphism $A \to U$ in $\Cc_X$, $j_{AU*}:\mathcal{S}(U) \to \mathcal{S}(A)$ is fully faithful.
\item A prestack $\mathcal{S}$ on $X$ is called a stack if for any $U \in \C_X$, and for any local isomorphism $A \to U$ in $\Cc_X$, $j_{AU*}:\mathcal{S}(U) \to \mathcal{S}(A)$ is an equivalence.
\end{itemize}
\end{df}

\begin{prop}\label{stack} Let $\mathcal{S}$ be a prestack on $X$. Then $\mathcal{S}$ is a
stack if and only if $\mathcal{S}$ satisfies the  following
conditions:
\begin{itemize}
\item[(i)] $\mathcal{S}$ is separated,
\item[(ii)] for any $U \in \C_X$ and for any local isomorphism $A\to U$  the restriction
functor $j_{AU*}:\mathcal{S}(U)\to\mathcal{S}(A)$ admits a left
adjoint $\imin {j_{AU}}$ satisfying $j_{AU*}\circ\imin
{j_{AU}}\simeq\id$ (or, equivalently, the functor $\imin {j_{AU}}$
is fully faithful).
\end{itemize}
\end{prop}
\dim\ \ The result follows from the fact that two categories are
equivalent if and only if they admit a pair of fully faithful
adjoint functors.

%(a) Assume that $\mathcal{S}$ satisfies (i) and (ii). Let $A\to V$ be a
%local isomorphism. We have a pair of adjoint functors
%\begin{equation*}
%\xymatrix{\mathcal{S}(A)   \ar@ <2pt> [r]^{\imin{j_{AV}}} &
%  \mathcal{S}(V) \ar@ <2pt> [l]^{j_{AV*}}. }
%\end{equation*}

%By (i) $j_{AV*}$ is fully faithful, and by (ii) $\imin {j_{AV}}$
%is fully faithful too. This implies that $j_{AV*}$ and $\imin
%{j_{AV}}$ are equivalence of categories inverse to each others.\\

%By (ii) we have $j_{AV*}\circ\imin {j_{AV}}\simeq\id$. Let us show
%$\imin {j_{AV}}\circ j_{AV*}\simeq\id$. Let $F,G \in \mathcal{S}(V)$.
%Since $\mathcal{S}$ is separated we have the chain of isomorphisms
%\begin{eqnarray*}
%\Ho_{\mathcal{S}(V)}(F,G) & \simeq & \Gamma(V;\ho_{\mathcal{S}}(F,G)) \\
%& \simeq & \Gamma(A;\ho_{\mathcal{S}}(F,G)) \\
%& \simeq & \Ho_{\mathcal{S}(V)}(j_{AV*}F,j_{AV*}G) \\
%& \simeq & \Ho_{\mathcal{S}(V)}(\imin {j_{AV}}j_{AV*}F,G).
%\end{eqnarray*}

%(b) Assume that $\mathcal{S}$ is a stack. Let $A\to V$ be a local
%isomorphism. Then the functor $j_{AV*}:\mathcal{S}(V) \to \mathcal{S}(A)$ is an
%equivalence, hence it admits a quasi-inverse $\imin {j_{AV}}$
%satisfying (ii). \nopagebreak \newline
\qed

%\begin{oss} During the proof of Proposition \ref{stack} we also
%showed that if $\mathcal{S}$ is separated and if for any $V \in \C_X$ and
%for any local isomorphism $A\to V$  the restriction functor
%$j_{AV*}:\mathcal{S}(V)\to\mathcal{S}(A)$ admits a left adjoint $\imin {j_{AV}}$,
%then $\imin {j_{AV}} \circ j_{AV*} \simeq \id$.
%\end{oss}

\section{Proper stacks}

Let $\C_X$ be a small category.
% admitting finite fiber products.
 In
this section we extend a result of \cite{KS01} to the case of a
site $X$ associated to a small category $\C_X$.

Let $\mathcal{S}$ be a prestack on $X$ and assume the following
hypothesis
\begin{equation}\label{j-1UV}
\begin{cases}
%\begin{array}{l}
\text{- for any $U,V \in \C_X$ and any morphism ${U}\to{V}$ in
$\Cc_X$, the functor} \\ \text{ \, $j_{UV*}:\mathcal{S}(V) \to
\mathcal{S}(U)$ admits a left adjoint $\imin {j_{UV}}$ satisfying} \\
\text{ \, $\id_{\mathcal{S}(U)}\iso{j_{UV*} \circ \imin
{j_{UV}}}$ (or, equivalently, $\imin {j_{UV}}$ is fully faithful),}\\
\text{- for all $U \in \C_X$ the category $\mathcal{S}(U)$ admits
small inductive limits.}
%\end{array}
\end{cases}
\end{equation}

\begin{lem}\label{AV} Let $\mathcal{S}$ be a prestack and assume \eqref{j-1UV}. Let $A \in \C_X^\wedge$ and $A \to V$. Then the functor $j_{AV*}$
admits a left adjoint, denoted by $\imin {j_{AV}}$.
\end{lem}
\dim\ \ Let $F=\{F_U\}_{(U \to A) \in \C_A} \in \mathcal{S}(A)$,
and let $\imin {j_{AV}}F:=\lind {(U \to A) \in \C_A} \imin
{j_{UV}}F_U$. This defines a functor $\imin
{j_{AV}}:\mathcal{S}(A)\to\mathcal{S}(V)$. Let $G \in
\mathcal{S}(V)$. We have the chain of isomorphisms
\begin{eqnarray*}
\Ho_{\mathcal{S}(V)}(\imin {j_{AV}}F,G) & = &
\Ho_{\mathcal{S}(V)}(\lind {(U
\to A) \in \C_A} \imin {j_{UV}}F_U,G) \\
& \simeq & \lpro {(U
\to A) \in \C_A}\Ho_{\mathcal{S}(V)}( \imin {j_{UV}}F_U,G) \\
& \simeq & \lpro {(U
\to A) \in \C_A}\Ho_{\mathcal{S}(U)}(F_U,j_{UV*}G) \\
& \simeq & \Ho_{\mathcal{S}(A)}(F,j_{AV*}G)
\end{eqnarray*}
\nopagebreak \qed

\begin{lem}\label{errou}  Let $\mathcal{S}$ be a prestack on $X$ satisfying \eqref{j-1UV}, let $U',U,V \in \C_X$
and ${U'}\to{U}\to{V}$. Then
\begin{itemize}
\item[(i)] there exists a canonical morphism
$\imin {j_{U'V}} \circ j_{U'V*}  \to \imin {j_{UV}} \circ
j_{UV*},$
\item[(ii)] we have $\imin {j_{U'V}} \circ j_{U'V*}  \simeq \imin {j_{U'V}} \circ
j_{U'V*} \circ \imin {j_{UV}} \circ j_{UV*}.$
\end{itemize}
\end{lem}
\dim\ \ (i) The adjunction morphism $\imin {j_{U'U}} \circ
j_{U'U*} \to \id_{\mathcal{S}(U)}$
 defines
$$\imin {j_{U'V}} \circ j_{U'V*} \simeq \imin {j_{UV}} \circ \imin {j_{U'U}} \circ j_{U'U*} \circ j_{UV*} \to \imin {j_{UV}} \circ j_{UV*}.$$

(ii) We have $j_{U'V*} \simeq j_{U'U*} \circ j_{UV*}$, and then
$$j_{U'V*} \circ  \imin {j_{UV}} \simeq j_{U'U*} \circ j_{UV*} \circ  \imin {j_{UV}} \simeq
j_{U'U*}.$$

Hence we have the chain of isomorphisms
$$\imin {j_{U'V}} \circ j_{U'V*} \circ \imin {j_{UV}} \circ j_{UV*}
\simeq \imin {j_{U'V}} \circ j_{U'U*} \circ j_{UV*} \simeq \imin
{j_{U'V}} \circ j_{U'V*}.$$ \nopagebreak \newline \qed

\begin{lem}  Let $\mathcal{S}$ be a prestack on $X$ satisfying \eqref{j-1UV}. Let $U,V,W \in \C_X$ and let $U \to W$, $V\to W$ be morphisms. Consider the diagram
$$ \xymatrix{ U \times_W V \ar[d] \ar[r] & V \ar[d] \\ U \ar[r] & W}
$$
where $U \times_W V \in \Cc_X$. Then there exists a canonical
morphism
\begin{equation}\label{UcapV}
\imin {j_{U \times_WVW}} \circ j_{U \times_WVW*} \to \imin
{j_{UW}} \circ j_{UW*}\circ \imin {j_{VW}} \circ j_{VW*}.
\end{equation}
\end{lem}
\dim\ \
%By Corollary \ref{VUV} we have an isomorphism
%\begin{eqnarray*}
%\imin {j_{U \times_XVW}} \circ j_{U \times_XVW*} & \simeq & \lind
%{W' \to U \times_XV \in \C_{U \times_XV}} \imin {j_{W'W}} \circ j_{W'W*} \\
%& \iso & \lind {W' \to U \times_XV \in \C_{U \times_XV}} \imin
%{j_{W'W}} \circ j_{W'W*} \circ \imin {j_{VW}} \circ j_{VW*}.
%\end{eqnarray*}
Since $U \times_W V \in \Cc_X$ for each $F \in \mathcal{S}(W)$ we
have
$$j_{U \times_W VW*}F=\{j_{W'W*}F\}_{(W' \to U \times_W V) \in
\C_{U \times_WV}} \in \mathcal{S}(U \times_W V)$$ hence as in
Lemma \ref{AV}
$$\imin {j_{U
\times_W VW}}j_{U \times_W VW*}F \simeq \lind {(W' \to U \times_W
V) \in \C_{U \times_W V}}\imin {j_{W'W}}j_{W'W*}F.$$ By Lemma
\ref{errou} we have $\imin {j_{W'W}} \circ j_{W'W*}\circ \imin
{j_{VW}} \circ j_{VW*} \simeq \imin {j_{W'W}} \circ j_{W'W*}$ for
each $(W' \to U \times_W V) \in \C_{U \times_W V}$. We have
natural morphisms
\begin{eqnarray*}
\imin {j_{U \times_WVW}} \circ j_{U \times_WVW*} & \iso & \imin
{j_{U \times_WVW}} \circ j_{U \times_WVW*}\circ \imin {j_{VW}}
\circ j_{VW*}\\
& \to & \imin {j_{UW}} \circ j_{UW*}\circ \imin {j_{VW}} \circ
j_{VW*}.
\end{eqnarray*} \nopagebreak  \qed

Let $U,V,W \in \C_X$ and let $U \to W$, $V\to W$ be morphisms. The
morphism \eqref{UcapV} induces a natural arrow
\begin{equation*}
\imin {j_{U \times_W V V}} \circ j_{U \times_W V U*} \iso  j_{VW*}
\circ \imin {j_{U \times_W VW}} \circ j_{U \times_W VW*} \circ
\imin {j _{UW}} \to j_{VW*} \circ \imin {j_{UW}}.
\end{equation*}

\begin{df} A proper stack $\mathcal{S}$ on $X$ is a prestack satisfying

\begin{itemize}

\item[PRS1] $\mathcal{S}$ is separated,

\item[PRS2] for each ${U}\in\C_X$, $\mathcal{S}(U)$ admits small inductive limits,
% and small projective limits,

\item[PRS3] for all $U,V \in \C_X$ and ${U}\to{V}$ the functor $j_{UV*}:\mathcal{S}(V) \to \mathcal{S}(U)$  commutes with $\Lind$,

\item[PRS4] for all $U,V \in \C_X$ and ${U}\to{V}$ the functor $j_{UV*}:\mathcal{S}(V) \to \mathcal{S}(U)$ admits a left adjoint $\imin {j_{UV}}$, satisfying $\id_{\mathcal{S}(U)}\iso{j_{UV*} \circ \imin
{j_{UV}}}$ (or, equivalently, the functor $\imin {j_{UV}}$ is
fully faithful),
%Moreover $\imin {j_{UV}}$ is exact.

%\item[(iv)] for all $U,V \in \C_X$ with ${U}\to{V}$, the functor $i_{UV!}$ is exact.

%\item[(v)] for all $U_1,U_2,V \in \C_X$ with $U_i \to V$, $i=1,2$
%and $\Ho_{\C_X}(U_1,U_2)=\Ho_{\C_X}(U_2,U_1)=\varnothing$ we have
%$j_{U_1U_2*} \circ \imin {j_{U_1U_2}} \simeq 0$

\item[PRS5]  for all ${V},{U},{W}\in\C_X$, $U \to W$ and $V \to W$, the morphism
$$\imin {j_{U \times_W V V}} \circ j_{U \times_W V U*} \to j_{VW*} \circ \imin
{j_{UW}}$$ is an isomorphism.

\end{itemize}
\end{df}

\begin{oss} Here $U \times_W V \in \Cc_X$, since we have not
assumed that $\C_X$ admits fiber products.
\end{oss}

\begin{lem} Let
us consider the following diagram
$$ \xymatrix{ A \times_V U \ar[d] \ar[r] & A \ar[d] \\ U \ar[r] & V}
$$
where $U,V \in \C_X$ and $A \in \C_X^\wedge$.
%with $U \to V$  and $A \to V$.
Let $\mathcal{S}$ be a proper stack on $X$. Then we have
$$j_{UV*} \circ \imin {j_{AV}} \simeq \imin j_{A \times_VUU} \circ
j_{A \times_V UA*}.$$
\end{lem}
\dim\ \ Let $F=\{F_{W}\}_{(W \to A) \in \C_A} \in \mathcal{S}(A)$.
We have the chain of isomorphisms
\begin{eqnarray*}
j_{UV*} \circ \imin {j_{AV}}F & \simeq & j_{UV*} \lind {(W \to A)
\in
\C_A} \imin {j_{WV}}F_W \\
& \simeq & \lind {(W \to A) \in
\C_A} j_{UV*} \imin {j_{WV}}F_W \\
& \simeq & \lind {(W \to A) \in \C_A} \imin {j_{U \times_V WU}}
j_{U
\times_V WW*}F_W \\
& \simeq & \lind {(W \to A) \in \C_A} \lind {(W' \to W \times_V U)
\in \C_{W \times_V U}}\imin {j_{W'U}} F_{W'}  \\
& \simeq & \lind {(W'' \to A \times_V U) \in \C_{A \times_V U}}
\imin {j_{W''U}}F_{W''},
\end{eqnarray*}
where the second and the third isomorphism follow from PRS3 and
 PRS5 respectively. The fourth isomorphism follows since $W
 \times_V U \in \Cc_X$ and we have
$$j_{U
\times_V WW*}F_W \simeq \{j_{W'W*}F_W\}_{(W' \to W \times_V U) \in
\C_{W \times_V U}} \simeq \{F_{W'}\}_{(W' \to W \times_V U) \in
\C_{W \times_V U}}.$$ On the other hand we have $j_{A \times_V
UA*}F \simeq \{F_{W''}\}_{(W'' \to A \times_V U) \in \C_{A
\times_V U}}$, hence
$$\imin j_{A \times_VUU} \circ
j_{A \times_V UA*}F \simeq \lind {(W'' \to A \times_V U) \in \C_{A
\times_V U}} \imin {j_{W''U}}F_{W''}.$$ \nopagebreak  \qed

\begin{teo} Let $X$ be a site associated to a small category $\C_X$. Let $\mathcal{S}$ be a proper stack on $X$. Then $\mathcal{S}$ is a
stack.
\end{teo}
\dim\ \ Let $A \to V$ be a local isomorphism. By Proposition
\ref{stack} it is enough to show that $j_{AV*}\circ\imin
{j_{AV}}\simeq\id$. Let $F=\{F_{V_i}\}_{(V_i \to A) \in \C_A} \in
\mathcal{S}(A)$. It satisfies, for each $V_i \to V_j$
\begin{equation} \label{isoij}
j_{V_iV_j*}F_{V_j} \iso F_{V_i}.
\end{equation}
We have to show that $j_{V_iV*}\imin {j_{AV}}F \simeq F_{V_i}$ for
each $V_i \to A$. Let us consider $V_{i_0} \to A$. By PRS5 and
\eqref{isoij}, for each $V_k \to A$ we have the chain of
isomorphisms
$$j_{V_{i_0}V*}\imin {j_{V_kV}}F_{V_k} \simeq \imin {j_{V_k
\times_V V_{i_0}V_{i_0}}}j_{V_k \times_V V_{i_0}V_k*}F_{V_k}
\simeq \imin {j_{V_k \times_V V_{i_0}V_{i_0}}}j_{V_k \times_V
V_{i_0}V_{i_0}*}F_{V_{i_0}}.$$ Hence we obtain the isomorphism
$$j_{V_{i_0}V*}\imin {j_{AV}}F \simeq \imin {j_{A \times_V
V_{i_0}V_{i_0}}}j_{A \times_V V_{i_0}V_{i_0}*}F_{V_{i_0}},$$ and
$\imin {j_{A \times_V V_{i_0}V_{i_0}}}j_{A \times_V
V_{i_0}V_{i_0}*}F_{V_{i_0}} \simeq F_{V_{i_0}}$ since
$\mathcal{S}$ is separated and $A \times_V V_{i_0} \to V_{i_0}$ is
a local isomorphism. \nopagebreak \newline \qed

\begin{es} Let $k$
be a field, and $X$ a topological space (or, more generally, let
$X$ be a site associated to an ordered-set category). The prestack
associating to an open set $U$ of $X$ the category of sheaves of
$k$-vector spaces\footnote[2]{More generally, one can consider
sheaves with values in a category $\A$ admitting small inductive
and projective limits, such that filtrant inductive limits are
exact and satisfying the ICP property (see \cite{KS} for a
detailed exposition).} on $U$ is a proper stack.
\end{es}

\end{document}